\documentclass[letterpaper]{article}
\usepackage{proceed2e}
\usepackage[margin=1in]{geometry}

\usepackage{times}

\usepackage{amsmath,amsthm,amssymb,mathtools}
\usepackage{tabu}
\usepackage{multirow}
\usepackage{url}
\usepackage{stackengine}
\usepackage{tikz}
\usetikzlibrary{arrows.meta,calc,decorations.pathmorphing}

\usepackage[round]{natbib} 
\DeclareRobustCommand{\VAN}[3]{#2 #1}

\renewcommand{\eqref}[1]{(\ref{eq:#1})}

\newtheorem{theorem}{Theorem}
\newtheorem{corollary}[theorem]{Corollary}

\newtheorem{proposition}[theorem]{Proposition}

\theoremstyle{definition} 

\newtheoremstyle{contstyle}
  {\topsep}   
  {\topsep}   
  {\normalfont} 
  {0pt}       
  {\bfseries} 
  {.}         
  {5pt plus 1pt minus 1pt} 
  {\thmname{#1}\thmnote{ #3} (continued)} 
\theoremstyle{contstyle}

\DeclarePairedDelimiterX\set[1]\lbrace\rbrace{\def\given{\;\delimsize\vert\;}#1}

\newcommand{\myempty}{\varnothing} 
\newcommand{\isdef}{\vcentcolon\protect\nolinebreak\mkern-1.2mu=}

\newcommand{\R}{\mathbb{R}} 

\DeclareMathOperator{\Cov}{Cov}

\newcommand{\trans}{^T}
\newcommand{\invtrans}{^{-T}}

\newcommand{\pa}{\mathrm{pa}}
\newcommand{\sib}{\mathrm{sib}}
\newcommand{\htr}{\mathrm{htr}}
\newcommand{\PD}{\mathrm{PD}}
\newcommand{\bA}{\mathbf{A}}
\newcommand{\bb}{\mathbf{b}}

\newcommand{\cI}{\ensuremath{\mathcal I}}
\newcommand{\cM}{\ensuremath{\mathcal M}}

\newcommand{\cY}{\ensuremath{\mathcal Y}}

\newcommand{\eqly}{algebraically} 
\newcommand{\eqic}{algebraic} 

\newcommand{\Eqic}{Algebraic}

\newcommand{\EQIC}{ALGEBRAIC}

\title{\Eqic{} Equivalence of Linear Structural Equation Models}
\author{ {\bf Thijs van Ommen} \\
Informatics Institute \\
University of Amsterdam\\
The Netherlands \\
\texttt{t.vanommen@uva.nl}
\And
{\bf Joris M. Mooij} \\
Informatics Institute \\
University of Amsterdam\\
The Netherlands \\
\texttt{j.m.mooij@uva.nl}
}
\begin{document}

\maketitle

\begin{abstract}
  Despite their popularity, many questions about the algebraic constraints imposed by linear structural equation models remain open problems. For causal discovery, two of these problems are especially important: the enumeration of the constraints imposed by a model, and deciding whether two graphs define the same statistical model. We show how the half-trek criterion can be used to make progress in both of these problems. We apply our theoretical results to a small-scale model selection problem, and find that taking the additional algebraic constraints into account may lead to significant improvements in model selection accuracy. 
\end{abstract}

\section{INTRODUCTION}





\begin{figure}[t]
  \centering
  \begin{tikzpicture}
    \node [circle,fill=black,inner sep=1pt] (a) at (0.866025403784,1.0) [label=90.0:$\mathstrut a$] {};
    \node [circle,fill=black,inner sep=1pt] (b) at (0.866025403784,0.0) [label=270.0:$\mathstrut b$] {};
    \node [circle,fill=black,inner sep=1pt] (c) at (0.0,0.5) [label=180.0:$\mathstrut c$] {};
    \node [circle,fill=black,inner sep=1pt] (d) at (1.73205080757,0.5) [label=0.0:$\mathstrut d$] {};
    \draw [blue,arrows=
      {_-Stealth[sep,length=1ex]}]
      (a) -- (b);
    \draw [red,dashed,arrows=
      {Stealth[sep,length=1ex]-Stealth[sep,length=1ex]}]
      (a) -- (c);
    \draw [red,dashed,arrows=
      {Stealth[sep,length=1ex]-Stealth[sep,length=1ex]}]
      (a) -- (d);
    \draw [red,dashed,arrows=
      {Stealth[sep,length=1ex]-Stealth[sep,length=1ex]}]
      (b) -- (c);
    \draw [blue,arrows=
      {_-Stealth[sep,length=1ex]}]
      (b) -- (d);
  \end{tikzpicture}
\caption{A mixed graph imposing the non-independence constraint \eqref{strangeconstraint} on the observational distribution.}\label{fig:strangeconstraint}
\end{figure}
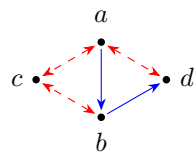

In a linear structural equation model (SEM), each 
variable of interest is a linear function of the other variables and a noise term, with possibly correlated noise terms.
Linear SEMs are popular in many fields of science, in no small part due to their causal interpretability \citep{SpirtesGlymourScheines2000,Pearl2000}. 
However, many questions about these models remain unanswered. For example, it is known that many of these models impose equality constraints on the observational distribution which do not correspond to (conditional) independences \citep{RichardsonSpirtes2002_MAGs}.
One example of these is the Verma constraint \citep{Robins1986_VermaConstraint,VermaPearl1991_VermaConstraint}.
But no general method exists which enumerates all constraints that hold in a model given its graphical representation. 
Figure~\ref{fig:strangeconstraint} shows another example of a graph that imposes an equality constraint in the linear case; 
for this type of constraint, no systematic approach exists yet.

Relatedly, given graphical representations of two models, it is often unclear whether the models can be distinguished based on observational data alone.
Both of these problems are great impediments to the development of methods that learn the structure of a linear SEM from observational data: constraint-based methods (e.g.~PC and FCI \citep{SpirtesGlymourScheines2000}) 
cannot test for constraints that are not yet well understood and may thus miss signals in the data, while score-based methods (e.g.~\citep{Chickering2002_GES,SilvaGhahramani2006}) would currently require the scoring of many models that then turn out to be indistinguishable.

The theoretical results of this paper achieve progress in addressing both problems mentioned above. For example, we provide an efficient method to find the constraint imposed by the graph in Figure~\ref{fig:strangeconstraint}, as well as similar constraints for many other graphs. Our results apply to models with latent confounders (represented by bidirected edges), including confounders between nodes that are also related by a direct causal effect (a structure known as a \emph{bow}), and even to models with directed cycles. We show how these results enable practical improvements on
model selection problems.

Suppose we see that our observational data obeys a non-independence equality constraint, such as the one imposed by the graph in Figure~\ref{fig:strangeconstraint}, and no other constraints. Then we can often draw very specific conclusions about the graph structure. Without knowledge of these constraints, automated methods for causal discovery would likely select a \emph{saturated} model 
(one imposing no constraints). This tells us very little about the graph structure, so e.g.~it does not allow us to predict the results of interventions. Constraint-based methods generally use only (conditional) independence constraints, so they would not be able to draw any conclusions in the situation described here. For purposes of model selection, we are thus interested in a notion of model equivalence that is more fine-grained than Markov equivalence (which only takes conditional independence constraints into account), yet not so fine-grained as to be impractical

The equivalence concept we propose in this paper is \emph{\eqic{} equivalence}: Two linear structural equation models are \eqly{} equivalent if they impose the same algebraic (i.e.~equality) constraints on the observational distribution.%
%
\footnote{Using terminology from algebraic geometry, the statistical models have the same Zariski closure \citep{CoxLittleOShea2015}.} %
These constraints take the form of polynomial equations over covariances $\sigma_{vw}$ of the observed variables.
Because the graph in Figure~\ref{fig:strangeconstraint} imposes such a constraint, it is not \eqly{} equivalent to the saturated model on four nodes, so a model selection method based on \eqic{} equivalence is able to distinguish the two, while a method based on Markov equivalence is not.

\begin{figure}[t]
  \centering
  \stackunder{\begin{tikzpicture} 
    \node [circle,fill=black,inner sep=1pt] (a) at (0,0) [label=270.0:$\mathstrut a$] {};
    \node [circle,fill=black,inner sep=1pt] (b) at (1,0) [label=270.0:$\mathstrut b$] {};
    \node [circle,fill=black,inner sep=1pt] (c) at (2,0) [label=270.0:$\mathstrut c$] {};
    \draw [blue,arrows=
      {Stealth[sep,length=1ex]-_}]
      (c) -- (b);
    \draw [red,dashed,arrows=
      {Stealth[sep,length=1ex]-Stealth[sep,length=1ex]}]
      (b) .. controls +(.25, .3) and +(-.25, .3) .. (c);
    \draw [blue,arrows=
      {Stealth[sep,length=1ex]-_}]
      (b) -- (a);
  \end{tikzpicture}}{(a)}
  \qquad
  \stackunder{\begin{tikzpicture} 
    \node [circle,fill=black,inner sep=1pt] (a) at (0,0) [label=270.0:$\mathstrut a$] {};
    \node [circle,fill=black,inner sep=1pt] (b) at (1,0) [label=270.0:$\mathstrut b$] {};
    \node [circle,fill=black,inner sep=1pt] (c) at (2,0) [label=270.0:$\mathstrut c$] {};
    \draw [blue,arrows=
      {Stealth[sep,length=1ex]-_}]
      (c) -- (b);
    \draw [blue,arrows=
      {Stealth[sep,length=1ex]-_}]
      (b) -- (a);
    \draw [blue,arrows=
      {_-Stealth[sep,length=1ex]}]
      (a) .. controls +(.5, .3) and +(-.5, .3) .. (c);
  \end{tikzpicture}}{(b)}
\caption{Two graphs whose statistical models are almost, but not entirely identical: (a) the instrumental variable model; (b) a saturated model.}\label{fig:instrumental}
\end{figure}
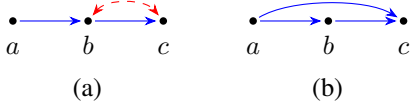

For an example of models not distinguished by \eqic{} equivalence, consider the mixed graph in Figure~\ref{fig:instrumental}(a), often called the instrumental variable model. 
This model contains all multivariate Gaussian distributions on the three variables
with $\sigma_{ab} \neq 0$, 
but excludes some with $\sigma_{ab} = 0$.
Because it imposes no equality constraints, it is \eqly{} equivalent to the saturated model on three nodes, represented for example by the graph in Figure~\ref{fig:instrumental}(b). The difference between the two models is a measure zero subset of their union, so that in a model selection problem, it would rarely be possible to distinguish between these models based on observational data alone. Thus it is appropriate that our proposed equivalence concept treats these models as equivalent.


By considering only equality constraints, we are also treating models as equivalent if they differ only by inequality constraints (of the form $f(\Sigma) \geq 0$).
%
It is known that linear Gaussian models with latent variables may impose inequality constraints. For example, consider the graph
with
three observed variables and a latent confounder between every pair of observed variables. 
The corresponding mixed graph has three nodes and a bidirected edge between each pair of nodes. This corresponds to the saturated model.
However, if we consider the entire graph on 6 nodes and then marginalize out the latent variables, the resulting set of distributions obeys a nontrivial inequality constraint \citep{RichardsonSpirtes2002_MAGs}. By using mixed graphs instead of including latent variables in our models explicitly, we are already simplifying away such inequality constraints. However, as we will show in Section~\ref{sec:ineq}, using mixed graphs does not get rid of all inequality constraints. We found that for models imposing an inequality constraint, maximum likelihood estimation can be challenging, so score-based methods for causal discovery would benefit significantly from being able to ignore these models. \Eqic{} equivalence provides a way to do this, as many models that impose an inequality constraint are \eqly{} equivalent to a model that imposes no such constraints.\footnote{In a model selection task, we can of course still check for inequality constraints, but we propose to do this \emph{after} selecting the \eqic{} equivalence class.}

\subsection{RELATIONS TO OTHER TYPES OF CONSTRAINTS}

When defining linear structural equation models statistically, often the noise terms are chosen to be Gaussian. For Gaussian variables, independence is equivalent to having zero covariance (and in turn to having zero correlation), and conditional independence is equivalent to zero \emph{partial correlation}. Because we do not want to assume that the data are generated by a Gaussian distribution, we need to distinguish between conditional independence and vanishing partial correlation. The equality constraints we consider express vanishing partial correlations, not conditional independences, so we will use that terminology from now on.

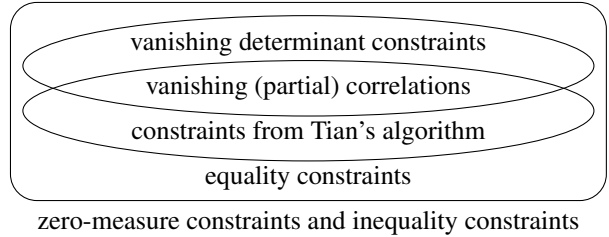
\begin{figure}[t]
  \centering
  \begin{tikzpicture} 
    \node (step) at (0, 17pt) {};
    \draw [rounded corners=10pt] ($(-112.7pt,0)-2*(step)$) rectangle (112.7pt,41pt);
    \draw (step) ellipse [x radius=108pt, y radius=19pt];
    \draw (0,0) ellipse [x radius=108pt, y radius=19pt];
    \node at ($1.5*(step)$) {vanishing determinant constraints};
    \node at ($0.5*(step)$) {vanishing (partial) correlations};
    \node at ($-0.5*(step)$) {constraints from Tian's algorithm};
    \node at ($-1.5*(step)$) {equality constraints};
    \node at ($-2.5*(step)$) {zero-measure constraints and inequality constraints};
  \end{tikzpicture}
\caption{The types of constraints imposed by linear structural equation models.}\label{fig:constraint_types}
\end{figure}
A partial correlation $\rho_{vw.S}$ is zero precisely when a certain submatrix of the observed covariance matrix has zero determinant. \citet{SullivantTalaskaDraisma2010} give a graphical characterization for the more general type of constraint where an arbitrary square submatrix has zero determinant; this also includes the vanishing tetrad constraints \citep{SpirtesGlymourScheines2000}. Together these are called \emph{vanishing determinant constraints} (see Figure~\ref{fig:constraint_types}).

The questions of constraint enumeration and model equivalence that we study here for the linear case, are studied for the general nonparametric case by \citet{TianPearl2002_TestableImplications} and \citet{ShpitserERR2014_NestedMarkovIntroduction}.
Tian's algorithm gives a sound enumeration of constraints in the general case, and 
\citet{Evans2015_arXiv_MarginsDiscrete} shows it to be complete in the discrete case. If other parametric assumptions are made, there may be additional constraints. Indeed this happens in the linear Gaussian case: Tian's algorithm returns no constraints for the graph in Figure~\ref{fig:strangeconstraint}, even though one exists.
%

Lists of algebraic constraints can also be obtained using algorithms from computer algebra, but these are in general impractically slow, sometimes taking many days even for very small graphs \citep{GarciapSpielvogelSullivant2010}. The methods we propose are based on the graphical criteria proposed by \citet{FoygelDraismaDrton2012_htc}, which can be checked in polynomial time.

The rest of this paper is structured as follows. Section~\ref{sec:prelim} discusses preliminaries about linear SEMs, the half-trek criterion, algebraic geometry, and our notation for sets of graphs. Our main theoretical contributions are in Section~\ref{sec:theoretical}, addressing the enumeration of algebraic constraints in Section~\ref{sec:allconstraints}, and a sufficient graphical criterion for \eqic{} equivalence in Section~\ref{sec:equiv}; further, Section~\ref{sec:ineq} gives an example of how inequality constraints may arise in linear SEMs. Experimental results demonstrating the practical usefulness of our results are presented in Section~\ref{sec:experimental}. Section~\ref{sec:conclusion} concludes the paper.
All proofs are in the supplementary material in Appendix~A, 
and a complete description of all \eqic{} equivalence classes on four nodes (acyclic) can be found in Appendix~B. 

\section{PRELIMINARIES}\label{sec:prelim}


We follow \citet{FoygelDraismaDrton2012_htc} for most of the notation defined in this section.

A \emph{mixed graph} $G = (V, D, B)$ consists of a set of nodes $V$, a set of directed edges $D$ which are ordered pairs of distinct nodes, and a set of bidirected edges $B$ which are unordered pairs of nodes. In this article, the word `graph' without qualification refers to mixed graphs. A node $x \in V$ with $(x,v) \in D$ is called a \emph{parent} of $v$, and the set of parents is denoted by $\pa(v)$; similarly, a node $x$ with $\set{x,v} \in B$ is called a \emph{sibling} of $v$ and the set is denoted $\sib(v)$. $G$ is called \emph{acyclic} if it contains no directed cycle (such a $G$ is also called \emph{acyclic directed mixed graph (ADMG)}). If $G$ is acyclic and contains no bidirected edges, it is called a \emph{directed acyclic graph (DAG)}.

Together with parameter vector $\lambda_0$ and parameter matrices $\Lambda$ and $\Omega$, the graph $G$ describes a distribution on observed variables $X$ via
\begin{equation*}
  X_v = \lambda_{0v} + \sum_{\mathclap{w \in \pa(v)}} \lambda_{wv} X_w + \epsilon_v
  \qquad \text{for $v \in V$,}
\end{equation*}
where the noise terms have covariances $\Cov(\epsilon_v, \epsilon_w) = \omega_{vw}$.
%
The parameter space is defined as follows.
Let $n = \lvert V \rvert$. $\R^D$ is the set of all $n \times n$ matrices $\Lambda$ with $\Lambda_{vw} \neq 0$ only if $(v,w) \in D$, and $\R^D_{\text{reg}}$ is the subset of $\R^D$ for which $I - \Lambda$ is invertible (for acyclic $G$, $\R^D = \R^D_{\text{reg}}$). $\PD_n$ is the set of all positive definite $n \times n$ matrices, and $\PD(B)$ is the subset consisting of all $\Omega$ with $\Omega_{vw} \neq 0$ only if $v = w$ or $\set{v,w} \in B$. The parameterization map $\phi_G$ maps parameters $(\Lambda, \Omega) \in \R^D_{\text{reg}} \times \PD(B)$ to covariance matrices $\Sigma \in \PD_n$ on the observed variables $X$ as follows:
\begin{equation*}
  \phi_G(\Lambda, \Omega) = (I - \Lambda)\invtrans \Omega (I - \Lambda)^{-1}.
\end{equation*}
The \emph{model} defined by a graph $G$ then consists of all covariance matrices $\Sigma$ that can be obtained for some setting of the parameters:
\begin{equation}\label{eq:model}
  \cM(G) \isdef \set{\phi_G(\Lambda, \Omega) \given (\Lambda, \Omega) \in \R^D_{\text{reg}} \times \PD(B)}.
\end{equation}
Note that the mean of $X$ can be set arbitrarily by choosing appropriate values for the parameter vector $\lambda_{0} \in \R^n$, regardless of the structure of $G$. Thus these aspects of the model carry no information for model selection, and we will ignore them here.

\subsection{THE HALF-TREK CRITERION}\label{sec:prelim_htc}

A central question about a mixed graph $G$ is that of \emph{(parameter) identifiability}: can the parameters $(\Lambda, \Omega)$ be uniquely recovered from $\Sigma$? A graph is called \emph{generically identifiable} (or almost-everywhere identifiable) if this is true of $\phi_G(\Lambda, \Omega)$ for all but a measure zero subset of $\R^D_{\text{reg}} \times \PD(B)$. Similarly, $G$ is called \emph{generically finite-to-one} if for almost all $(\Lambda, \Omega)$, the number of parameter values mapped by $\phi_G$ to the same $\Sigma$ is finite, and \emph{generically infinite-to-one} if this number is infinite for almost all parameter values. We will sometimes omit the qualifier `generically' when talking about (in)finite-to-one models.


\citet{FoygelDraismaDrton2012_htc} present two graphical criteria to decide in which of the above categories a graph $G$ belongs. A graph is called \emph{HTC-identifiable} if it meets the condition for being generically identifiable; \emph{HTC-nonidentifiable} if it meets the condition for being generically infinite-to-one; and \emph{HTC-inconclusive} otherwise. Because neither criterion is necessary, the class of HTC-inconclusive graphs contains generically identifiable, finite-to-one, and infinite-to-one graphs. While the criteria are not complete, they are quite powerful. For example, all bow-free acyclic graphs are HTC-identifiable (thus implying the earlier identification result of \citet{BritoPearl2002_SEM}), and so are many graphs containing bows or directed cycles.

The proof of HTC-identifiability in \citep{FoygelDraismaDrton2012_htc} is constructive: it gives an algorithm that, given $\Sigma \in \cM(G)$, computes parameters such that $\Sigma = \phi_G(\Lambda, \Omega)$ (except for a measure zero subset). To apply our Theorem~\ref{thm:allconstraints} below, some details of HTC-identifiability and this algorithm are needed; for the rest, we refer to \citep[proof of Theorem~1]{FoygelDraismaDrton2012_htc}.
A \emph{half-trek} from $v$ to $w$ is either a directed path, or a path consisting of one bidirected edge followed by directed edges towards $w$. We write $\htr(v)$ (\emph{half-trek reachable}) for the set of nodes that are reachable from $v$ by half-treks.\footnote{Here we follow the (more natural) definition of \citet{ChenTianPearl2014_LSEMOveridentify} rather than that of \citet{FoygelDraismaDrton2012_htc}.} HTC-identifiability of a graph $G$ requires that for each node $v \in V$, a set $Y_v \subseteq V \setminus (\set{v} \cup \sib(v))$ exists, consisting of nodes $y$ with $v \in \htr(y)$, and the set as a whole satisfying $\lvert Y_v \rvert = \lvert \pa(v) \rvert$.
There are some additional restrictions which we omit here, except to point out that for all $v, w \in V$, at most one of $v \in Y_w$ and $w \in Y_v$ can hold.
As an example, the instrumental variable model (Figure~\ref{fig:instrumental}(a)) is HTC-identifiable with $Y_a = \myempty$ and $Y_b = Y_c = \set{a}$.

Using these sets $Y_v$, the algorithm for finding $\Lambda$ solves a sequence of linear systems, one for each $v \in V$: Let $Y_v = \set{y_1, \ldots, y_n}$ and $\pa(v) = \set{p_1, \ldots, p_n}$, and define $\bA \in \R^{n \times n}$ and $\bb \in \R^n$ as
\begin{align*}
  \bA_{ij} &=
  \begin{cases*}
    [(I - \Lambda)\trans \Sigma]_{y_i p_j} & if $y_i \in \htr(v)$,\\
    \Sigma_{y_i p_j} & if $y_i \notin \htr(v)$;
  \end{cases*}\\
  \bb_{i} &=
  \begin{cases*}
    [(I - \Lambda)\trans \Sigma]_{y_i v} & if $y_i \in \htr(v)$,\\
    \Sigma_{y_i v} & if $y_i \notin \htr(v)$.
  \end{cases*}
\end{align*}
Then the vector $\Lambda_{\pa(v), v}$ is found by solving $\bA \cdot \Lambda_{\pa(v), v} = \bb$. After the entire matrix $\Lambda$ has been found this way, $\Omega$ is given by $(I - \Lambda)\trans \Sigma (I - \Lambda)$.

\subsection{SOME TERMINOLOGY FROM ALGEBRAIC GEOMETRY}\label{sec:prelim_alg_geom}

This section very briefly highlights the main terms from algebraic geometry; we refer to \citet{CoxLittleOShea2015} for further reading.

Algebraic geometry studies sets of points defined by systems of polynomial equations. For points in $\R^n$, such a set is called an \emph{affine variety}. The set of polynomials that are identically zero on some set of points form an \emph{ideal}. The \emph{Zariski closure} of a set of points is obtained by first finding the ideal of the set, then taking all points for which the polynomials in the ideal all vanish.


\subsection{NOTATION FOR SETS OF GRAPHS}\label{sec:graph_notation}

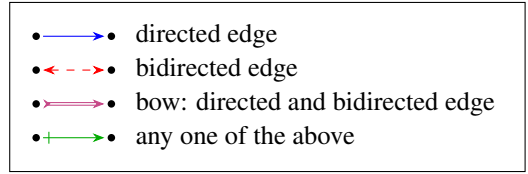
\begin{figure}[t]
  \centering
  \begin{tikzpicture} 
    \node (step) at (0, -.45) {};
    \draw ($(-10pt,0)-(step)$) rectangle ($(6.5,0)+4*(step)$);
    \node [circle,fill=black,inner sep=1pt] (a1) at (0,0) {};
    \node [circle,fill=black,inner sep=1pt] (b1) at (1,0) {};
    \node [circle,fill=black,inner sep=1pt] (a2) at ($(0,0)+(step)$) {};
    \node [circle,fill=black,inner sep=1pt] (b2) at ($(1,0)+(step)$) {};
    \node [circle,fill=black,inner sep=1pt] (a3) at ($(0,0)+2*(step)$) {};
    \node [circle,fill=black,inner sep=1pt] (b3) at ($(1,0)+2*(step)$) {};
    \node [circle,fill=black,inner sep=1pt] (a4) at ($(0,0)+3*(step)$) {};
    \node [circle,fill=black,inner sep=1pt] (b4) at ($(1,0)+3*(step)$) {};
    \draw [blue,arrows=
      {Stealth[sep,length=1ex]-_}]
      (b1) -- (a1);
    \node [anchor=west] at (1.2,0) {directed edge};
    \draw [red,dashed,arrows=
      {Stealth[sep,length=1ex]-Stealth[sep,length=1ex]}]
      (a2) -- (b2);
    \node [anchor=west] at ($(1.2,0)+(step)$) {bidirected edge};
    \draw [magenta!80!black,double,arrows=
      {Stealth[reversed,sep,length=.64ex,inset=.25ex,width=.8ex]-Stealth[sep,length=1ex]}]
      (a3) -- (b3);
    \node [anchor=west] at ($(1.2,0)+2*(step)$) {bow: directed and bidirected edge};
    \draw [green!67!black,arrows=
      {Stealth[sep,length=1ex]-Bar[sep,width=1ex] __}]
      (b4) -- (a4);
    \node [anchor=west] at ($(1.2,0)+3*(step)$) {any one of the above};
  \end{tikzpicture}
\caption{Legend for the edges we use to denote (sets of) graphs.}\label{fig:legend_short}
\end{figure}

The notation described here will be used in the remaining figures in this paper, and extensively in Appendix~B 
in the supplementary material. When drawing a mixed graph, we draw directed edges in solid blue and bidirected edges in dashed red (see Figure~\ref{fig:legend_short}; note that the colours are redundant but may aid visual distinction). Because we often want to show a set of graphs with a common node set $V$ (e.g., an \eqic{} equivalence class of graphs), we use some new notation to avoid listing all graphs one by one. This notation is based on the \emph{skeletons} of the graphs. Formally, the skeleton $S(G)$ of $G = (V,D,B)$ is the undirected graph on V that has an edge between a pair of nodes if there is at least one edge of any type between them in $G$. Similarly, our \emph{graph patterns} also have at most one edge between each pair of nodes, with different markings
to show what (combinations of) edges may occur between those nodes. For example, in this notation, a bow (a directed and a bidirected edge between two nodes) is represented visually as a double magenta line with `fletching' at the back.

In a pattern representing a set of graphs, other edges occur in places where the graphs differ.
The only such edge that appears in the patterns in the main paper is a green arrow with a plus sign as its tail. In all of these patterns, it can be understood to mean either a directed edge, a bidirected edge, or both.
Its meaning changes if another edge is incident at the endpoint with the plus sign; for this case, additional markings are used and defined in Appendix~B.

\section{THEORETICAL RESULTS}\label{sec:theoretical}

We now present our main theoretical results.

\subsection{ALGEBRAIC CONSTRAINTS}\label{sec:allconstraints}

The statistical model $\cM(G)$ for a graph $G$ is defined by \eqref{model} in terms of the parameterization $\phi_G$. This makes it hard to decide whether a given $\Sigma$ belongs to a model, or whether two models differ. For these purposes, a more usable description of $\cM(G)$ would be as a set of constraints that must hold for any $\Sigma \in \cM(G)$: a set of functions $f_1,\ldots,f_k: \Sigma \to \R$ such that
\begin{equation*}
  \cM(G) \subseteq \set{\Sigma \in \PD_{\lvert V \rvert} \given f_1(\Sigma) = \ldots = f_k(\Sigma) = 0},
\end{equation*}
with the set difference a measure zero subset of the right-hand set.%
\footnote{Put more precisely, the difference is contained in an affine variety of lower dimension.} 
For a DAG $G$, this can be done by choosing each $f_i$ to be a partial correlation: $f_i(\Sigma) = \rho_{{v_i w_i}.S_i}(\Sigma)$ \citep{RichardsonSpirtes2002_MAGs}. 
Put differently, vanishing partial correlation constraints fully describe such models. 

In the more general case where $G$ is allowed to contain bidirected edges and directed cycles, vanishing partial correlation constraints are not always expressive enough to describe $\cM(G)$ (as we saw in the example of Figure~\ref{fig:strangeconstraint}). Below we present a method that allows us to write down a list of equality constraints that together describe $\cM(G)$, for any \emph{HTC-identifiable graph} $G$ \citep{FoygelDraismaDrton2012_htc}.

Any model $\cM(G)$ is a \emph{semi-algebraic} subset of $\PD_{\lvert V \rvert}$: it can be described by a set of polynomial (in)equalities \citep{SullivantTalaskaDraisma2010}. Thus we may use the terms `equality constraint' and `algebraic (i.e.~polynomial) constraint' interchangeably, and similarly for `inequality' and `semi-algebraic'. 

The question of finding constraints that describe $\cM(G)$ is related to the question of parameter identifiability. Identifiability of a graph $G$ means that for given $\Sigma$, \emph{at most} one pair of parameter values $(\Lambda, \Omega)$ exists for which $\phi_G(\Lambda, \Omega) = \Sigma$ (in other words, the \emph{fibre} $\phi_G^{-1}(\Sigma)$ has cardinality at most one). We ask for what $\Sigma$ (up to measure zero subsets as above) \emph{exactly one} such pair exists. The main tool used here is the algorithm for retrieving the parameters of an HTC-identifiable graph that appears in the proof of Theorem~1 of \citet{FoygelDraismaDrton2012_htc}, and that we sketched in Section~\ref{sec:prelim_htc}. For given $G$ and $\cY = (Y_v)_v$ (which we will call \emph{HTC-identifying sets}) satisfying the conditions of the HTC-identifiability theorem, this algorithm defines a rational function that maps $\Sigma$ to $\Lambda$ by solving a sequence of linear equation systems. We write this function as $\Lambda_\cY(\Sigma)$.\footnote{Note that this function does not depend on the ordering $\prec$ that also appears in the HTC-identifiability theorem.} Points $\Sigma$ for which the algorithm encounters a singular matrix are excluded from the domain of $\Lambda_\cY$.

Theorem~\ref{thm:allconstraints} below shows that a graph $G$ imposes the following constraints:
\begin{multline}\label{eq:allconstraints}
  [(I - \Lambda_\cY(\Sigma))\trans \Sigma (I - \Lambda_\cY(\Sigma))]_{vw} = 0%
  \\ 
  \text{for all $\set{v, w} \notin B$ with $v \notin Y_w$ and $w \notin Y_v$.}
\end{multline}
Each left-hand side in these equations is a rational function of $\Sigma$: a function of the form $p(\Sigma) / q(\Sigma)$, with $p$ and $q$ polynomials. Instead of these rational constraints, it will often be useful to consider the polynomial constraints $p(\Sigma) = 0$ obtained by multiplying out the denominator. For $\Sigma$ with $q(\Sigma) \neq 0$, $p(\Sigma) = 0$ iff $p(\Sigma)/q(\Sigma) = 0$, so the two forms of the constraints agree whenever both are defined (see the proof of Theorem~\ref{thm:allconstraints} for details). The polynomial constraints have the advantage that they are defined everywhere.

As an example, consider again the graph in Figure~\ref{fig:strangeconstraint}. For $\cY$ with $Y_b = \set{a}$, $Y_d = \set{c}$ and $Y_a = Y_c = \myempty$, \eqref{allconstraints} gives us one rational constraint (for $v = b$, $w = d$):
\begin{equation*}
  \begin{bmatrix}
    1 & -\Lambda_\cY(\Sigma)_{ab}
  \end{bmatrix}
  \begin{bmatrix}
    \sigma_{bd} & \sigma_{bb} \\
    \sigma_{ad} & \sigma_{ab}
  \end{bmatrix}
  \begin{bmatrix}
    1 \\
    -\Lambda_\cY(\Sigma)_{bd}
  \end{bmatrix} = 0,
\end{equation*}
where $\Lambda_\cY(\Sigma)_{ab} = \sigma_{ab} / \sigma_{aa}$ and $\Lambda_\cY(\Sigma)_{bd} = \sigma_{cd} / \sigma_{bc}$. Multiplying out the denominators, we obtain the polynomial constraint
\begin{equation}\label{eq:strangeconstraint}
  \sigma_{aa} \sigma_{bd} \sigma_{bc} - \sigma_{aa} \sigma_{bb} \sigma_{cd} - \sigma_{ab} \sigma_{ad} \sigma_{bc} + \sigma_{ab}^2 \sigma_{cd} = 0.
\end{equation}
%
%
%
\begin{theorem}\label{thm:allconstraints}
  For an HTC-identifiable graph $G$ with HTC-identifying sets $\cY = (Y_v)_v$ and generic $\Sigma$, $\Sigma \in \cM(G)$ iff $\Lambda_\cY(\Sigma) \in \R^D_{\text{reg}}$ and $\Sigma$ satisfies the rational constraints \eqref{allconstraints}. %
%
  A stronger statement holds in one direction: \emph{All} (not merely generic) $\Sigma \in \cM(G)$ satisfy the polynomial constraints described above.
\end{theorem}
This means that the constraints \eqref{allconstraints} define the model $\cM(G)$ up to a measure zero set that may satisfy all constraints, but still be missing from $\cM(G)$. In particular, it shows that the model of an HTC-identifiable graph imposes no semi-algebraic (i.e.~inequality) constraints. This is not true for general graphs, as we will see in Section~\ref{sec:ineq}.

The HT-overidentifying constraints from \citet{ChenTianPearl2014_LSEMOveridentify} are also based on (a version of) the half-trek criterion, but unlike \eqref{allconstraints}, they do not give a full description of the algebraic constraints imposed by a graph $G$: no constraint is found for pairs $\set{v,w}$ with $v \not\in \htr(w)$ and $w \not\in \htr(v)$. An example is given in Appendix~\ref{app:proofs}.

There may often be several ways of expressing $\Lambda$ in terms of $\Sigma$, using different HTC-identifying sets $\cY$. As a result, the rational constraints may look very different, though they may become the same when converted to polynomial form. For example, a different choice of $\cY$ for the graph in Figure~\ref{fig:strangeconstraint} would have led to a different rational expression than we found above. Similarly, sets of constraints may be found which are different when compared one by one, even in polynomial form, but which together still define the same model (in the terminology of algebraic geometry, they generate the same ideal).
For example, if $\sigma_{ab} = 0$, then $\sigma_{ac} = 0$ iff $\rho_{ac.b} = 0$, so that a model satisfying all these constraints can be described in two ways using two constraints, or redundantly using all three constraints.
So while the set of constraints \eqref{allconstraints} is not unique, Theorem~\ref{thm:allconstraints} shows that they are a complete description of $\cM(G)$: the model imposes no algebraic constraints beyond these.

Theorem~\ref{thm:allconstraints} only applies to HTC-identifiable graphs. This immediately excludes graphs which are not generically identifiable; we will revisit those in Sections~\ref{sec:equiv} and~\ref{sec:ineq}. However, it also excludes graphs that are generically identifiable but not HTC-identifiable.
For the case of acyclic graphs on four variables, HTC-identifiability is complete (in the sense that all generically identifiable graphs are also HTC-identifiable), but this is no longer true when either more nodes or cycles are allowed. Finding more general sufficient graphical criteria for generic identifiability is the topic of ongoing research \citep{ChenTianPearl2014_LSEMOveridentify,Chen2016_Overidentify_NIPS,DrtonWeihs2016_HTCidWithAncestorDecomposition,ChenKB2017_TR}. Because these criteria are extensions of HTC-identifiability, our Theorem~\ref{thm:allconstraints} might be extended to work with these criteria as well.


\subsection{\EQIC{} EQUIVALENCE AND INFINITE-TO-ONE GRAPHS}\label{sec:equiv}

We now turn to the second central problem we address in this paper, namely checking whether two graphs $G_1$ and $G_2$ are \eqly{} equivalent.
Theorem~\ref{thm:allconstraints} from the previous section in principle allows us to do this: we need to check that each equality constraint imposed by $G_1$ is implied by $G_2$'s equality constraints, and vice versa.
However, the general solution for such tasks from algebraic geometry (computing Gr\"obner bases \citep{CoxLittleOShea2015}) is computationally extremely expensive.
The theorem below gives a sufficient condition for \eqic{} equivalence that relies only on graphical criteria that can be checked efficiently.

Another limitation of Theorem~\ref{thm:allconstraints} is that it only applies to HTC-identifiable graphs, so it does not help us in finding equivalences involving graphs that are not HTC-identifiable. This issue is also addressed by the following theorem.


\begin{theorem}\label{thm:equiv_sufficient}
  If $G$ is generically infinite-to-one, $G'$ generically finite-to-one and obtained by deleting $k$ edges from $G$, and this $k$ is the smallest number for which such a $G'$ exists, then $G$ and $G'$ are \eqly{} equivalent.
  If further $G'$ imposes no inequality constraint, $\cM(G)$ and $\cM(G')$ are equal up to a measure zero subset.
\end{theorem}
A relation between forms of model equivalence and lack of parameter identifiability has been shown previously by \citet{BekkerMerckensWansbeek1994}.
The following two corollaries serve to illustrate the power of this theorem.
\begin{corollary}\label{cor:inf_equiv}
  All infinite-to-one graphs are \eqly{} equivalent to some finite-to-one graph.
\end{corollary}
This implies that it is not a limitation that Theorem~\ref{thm:allconstraints} does not apply to infinite-to-one graphs: for any such graph, these always exists a finite-to-one graph that we can consider instead. In particular, if this finite-to-one graph is HTC-identifiable, then this gives us a complete description of the infinite-to-one graph as a set of equality constraints. 
\begin{corollary}\label{cor:fin_equiv}
  If two generically finite-to-one graphs $G'_1, G'_2$ are each obtained by deleting an edge from a single generically infinite-to-one graph $G$, then $G'_1$ and $G'_2$ are \eqly{} equivalent to each other (and their models are equal up to a measure zero subset if they impose no inequality constraints).
\end{corollary}
Combined with graphical criteria for generic (in)finite-to-oneness, this gives a sufficient \emph{graphical} condition for model equivalence. As an example, take $G'_1$ and $G'_2$ to be the two graphs in Figure~\ref{fig:instrumental}, and $G$ to be their union (thus having three nodes and four edges). $G'_1$ and $G'_2$ are generically finite-to-one, while $G$ is generically infinite-to-one (it is HTC-nonidentifiable). Then Corollary~\ref{cor:fin_equiv} states that $G'_1$ and $G'_2$ are \eqly{} equivalent. In fact, both graphs are HTC-identifiable, so they impose no inequality constraints, and their models are thus equal up to measure zero subsets.


By repeatedly applying Corollary~\ref{cor:fin_equiv}, equivalence of many more pairs of models may be established. This is illustrated by
the following proposition (which is similar to \citep[Theorem~2]{NowzohourMEB2017_arXiv}, though there the stronger relation of distributional equivalence is shown).
\begin{proposition}\label{prop:nowzohour}
  If two bow-free acyclic graphs have the same skeleton, and any 2-edge path through three distinct nodes that is a collider in one graph is a collider in both, the graphs are \eqly{} equivalent.
\end{proposition}

\subsubsection{\Eqic{} Equivalence Classes on Four Nodes}\label{sec:modelclasses}

Using Theorems~\ref{thm:allconstraints} and~\ref{thm:equiv_sufficient} and information on the identifiability of the acyclic graphs on four nodes, we can determine the \eqic{} equivalence classes of these graphs. An explicit description of these classes is given in Appendix~B. 
Here we describe some details of how these results were derived.


Using HTC-identifiability and HTC-nonidentifiability, almost all acyclic graphs on four nodes can be classified as either generically identifiable or generically infinite-to-one.
With the additional information from \citep{FoygelDraismaDrton2012_htc} that the remaining graphs are generically finite-to-one (see Section~\ref{sec:ineq}), we can
apply Theorem~\ref{thm:equiv_sufficient} to all these graphs. Together with transitivity, this partitions the set of graphs into 419 subsets, which we will call \emph{clusters} here. Because the theorem only gives a sufficient condition for \eqic{} equivalence, this partition may be finer than the partition into \eqic{} equivalence classes. We still need to check if \eqic{} equivalences exist between different clusters.

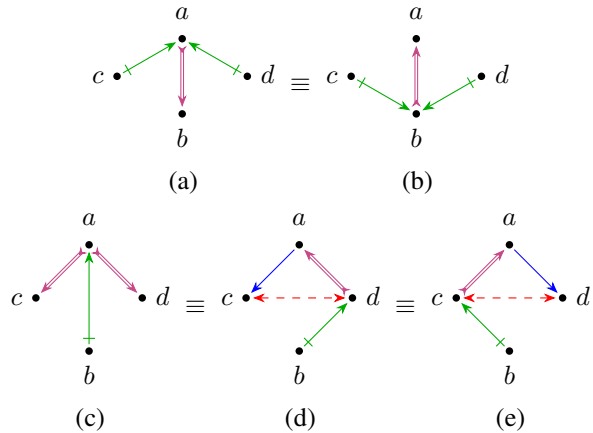
\begin{figure}[t]
  \centering
  \stackunder{\begin{tikzpicture}[baseline=.3cm] 
    \node [circle,fill=black,inner sep=1pt] (a) at (0.866025403784,1.0) [label=90.0:$\mathstrut a$] {};
    \node [circle,fill=black,inner sep=1pt] (b) at (0.866025403784,0.0) [label=270.0:$\mathstrut b$] {};
    \node [circle,fill=black,inner sep=1pt] (c) at (0.0,0.5) [label=180.0:$\mathstrut c$] {};
    \node [circle,fill=black,inner sep=1pt] (d) at (1.73205080757,0.5) [label=0.0:$\mathstrut d$] {};
    \draw [magenta!80!black,double,arrows=
      {Stealth[reversed,sep,length=.64ex,inset=.25ex,width=.8ex]-Stealth[sep,length=1ex]}]
      (a) -- (b);
    \draw [green!67!black,arrows=
      {Stealth[sep,length=1ex]-Bar[sep,width=1ex] __}]
      (a) -- (c);
    \draw [green!67!black,arrows=
      {Stealth[sep,length=1ex]-Bar[sep,width=1ex] __}]
      (a) -- (d);
  \end{tikzpicture}}{(a)}
  $\equiv$
  \stackunder{\begin{tikzpicture}[baseline=.3cm]
    \node [circle,fill=black,inner sep=1pt] (a) at (0.866025403784,1.0) [label=90.0:$\mathstrut a$] {};
    \node [circle,fill=black,inner sep=1pt] (b) at (0.866025403784,0.0) [label=270.0:$\mathstrut b$] {};
    \node [circle,fill=black,inner sep=1pt] (c) at (0.0,0.5) [label=180.0:$\mathstrut c$] {};
    \node [circle,fill=black,inner sep=1pt] (d) at (1.73205080757,0.5) [label=0.0:$\mathstrut d$] {};
    \draw [magenta!80!black,double,arrows=
      {Stealth[reversed,sep,length=.64ex,inset=.25ex,width=.8ex]-Stealth[sep,length=1ex]}]
      (b) -- (a);
    \draw [green!67!black,arrows=
      {Stealth[sep,length=1ex]-Bar[sep,width=1ex] __}]
      (b) -- (c);
    \draw [green!67!black,arrows=
      {Stealth[sep,length=1ex]-Bar[sep,width=1ex] __}]
      (b) -- (d);
  \end{tikzpicture}}{(b)}

  \stackunder{\begin{tikzpicture}[baseline=-.2cm]
    \node [circle,fill=black,inner sep=1pt] (a) at (0,.707) [label=90:$\mathstrut a$] {};
    \node [circle,fill=black,inner sep=1pt] (b) at (0,-.707) [label=-90:$\mathstrut b$] {};
    \node [circle,fill=black,inner sep=1pt] (c) at (-.707,0) [label=180:$\mathstrut c$] {};
    \node [circle,fill=black,inner sep=1pt] (d) at (.707,0) [label=0:$\mathstrut d$] {};
    \draw [green!67!black,arrows=
      {Stealth[sep,length=1ex]-Bar[sep,width=1ex] __}]
      (a) -- (b);
    \draw [magenta!80!black,double,arrows=
      {Stealth[reversed,sep,length=.64ex,inset=.25ex,width=.8ex]-Stealth[sep,length=1ex]}]
      (a) -- (c);
    \draw [magenta!80!black,double,arrows=
      {Stealth[reversed,sep,length=.64ex,inset=.25ex,width=.8ex]-Stealth[sep,length=1ex]}]
      (a) -- (d);
  \end{tikzpicture}}{(c)}
  $\equiv$
  \stackunder{\begin{tikzpicture}[baseline=-.2cm]
    \node [circle,fill=black,inner sep=1pt] (a) at (0,.707) [label=90:$\mathstrut a$] {};
    \node [circle,fill=black,inner sep=1pt] (b) at (0,-.707) [label=-90:$\mathstrut b$] {};
    \node [circle,fill=black,inner sep=1pt] (c) at (-.707,0) [label=180:$\mathstrut c$] {};
    \node [circle,fill=black,inner sep=1pt] (d) at (.707,0) [label=0:$\mathstrut d$] {};
    \draw [blue,arrows=
      {_-Stealth[sep,length=1ex]}]
      (a) -- (c);
    \draw [magenta!80!black,double,arrows=
      {Stealth[sep,length=1ex]-Stealth[reversed,sep,length=.64ex,inset=.25ex,width=.8ex]}]
      (a) -- (d);
    \draw [green!67!black,arrows=
      {__ Bar[sep,width=1ex]-Stealth[sep,length=1ex]}]
      (b) -- (d);
    \draw [red,dashed,arrows=
      {Stealth[sep,length=1ex]-Stealth[sep,length=1ex]}]
      (c) -- (d);
  \end{tikzpicture}}{(d)}
  $\equiv$
  \stackunder{\begin{tikzpicture}[baseline=-.2cm]
    \node [circle,fill=black,inner sep=1pt] (a) at (0,.707) [label=90:$\mathstrut a$] {};
    \node [circle,fill=black,inner sep=1pt] (b) at (0,-.707) [label=-90:$\mathstrut b$] {};
    \node [circle,fill=black,inner sep=1pt] (c) at (-.707,0) [label=180:$\mathstrut c$] {};
    \node [circle,fill=black,inner sep=1pt] (d) at (.707,0) [label=0:$\mathstrut d$] {};
    \draw [blue,arrows=
      {_-Stealth[sep,length=1ex]}]
      (a) -- (d);
    \draw [magenta!80!black,double,arrows=
      {Stealth[sep,length=1ex]-Stealth[reversed,sep,length=.64ex,inset=.25ex,width=.8ex]}]
      (a) -- (c);
    \draw [green!67!black,arrows=
      {__ Bar[sep,width=1ex]-Stealth[sep,length=1ex]}]
      (b) -- (c);
    \draw [red,dashed,arrows=
      {Stealth[sep,length=1ex]-Stealth[sep,length=1ex]}]
      (d) -- (c);
  \end{tikzpicture}}{(e)}
  \caption{The \eqic{} equivalences not detected by Theorem~\ref{thm:equiv_sufficient}: clusters (a) and~(b) are \eqly{} equivalent to each other, and the same is true for clusters (c), (d) and~(e).}\label{fig:extra_equiv}
\end{figure}

Among clusters imposing two or more equality constraints, all but three (up to graph isomorphism) can be described by vanishing (partial) correlation constraints. Because graphs in different Markov equivalence classes must also be in different \eqic{} equivalence classes, we only need to focus on these remaining three clusters. One of these imposes the Verma constraint \citep{VermaPearl1991_VermaConstraint}; the other two both impose the constraints $\sigma_{cd} = 0$ and $\sigma_{ac} \sigma_{bd} - \sigma_{ad} \sigma_{bc} = 0$ (a vanishing tetrad constraint). These two latter clusters are shown in Figure~\ref{fig:extra_equiv}(a) and~(b). So those two clusters are \eqly{} equivalent, while no other \eqic{} equivalences among these graphs were missed by Theorem~\ref{thm:equiv_sufficient}.

For graphs imposing only one equality constraint, it is much easier to check if two graphs impose the same equality constraint, as we do not need to worry about the possibility of two equality constraints implying a third (see Section~\ref{sec:allconstraints}). In algebraic terminology, these models are described by principal ideals; for these, equality can be checked by normalizing the generating polynomials so that their leading coefficients equal one \citep{CoxLittleOShea2015}. This way, we find that among graphs imposing one equality constraint, the three clusters (up to isomorphism) shown in Figure~\ref{fig:extra_equiv}(c), (d) and~(e) are actually \eqly{} equivalent to each other, while all others are different. This leaves a total of 389 \eqic{} equivalence classes.


Interestingly, the graphs imposing a vanishing tetrad constraint and no other equality constraints were all determined to be \eqly{} equivalent by Theorem~\ref{thm:equiv_sufficient}.


\subsubsection{Consequences for Model Selection}\label{sec:consequences}

The theoretical results we presented above offer enormous benefits to model selection, in particular to score-based methods. Without any knowledge of model equivalence, a score-based method for model selection may in principle need to score all different graphs. Even for the limited case of acyclic graphs on four nodes, there are 34752 distinct graphs, making such an approach clearly not practical. For this reason, most score-based methods for model selection rely on a coarser concept of equivalence such as Markov equivalence, or even limit themselves to DAGs, ignoring the possibility of latent confounders entirely \citep{DrtonMaathuis2016_review_StrucLearn}. As demonstrated in Section~\ref{sec:modelclasses}, our theoretical results can be used to determine the 389 \eqic{} equivalence classes. Scoring just one representative of each class leads to huge computational savings, making \eqic{} equivalence class selection feasible. As will be elaborated on in Section~\ref{sec:experimental}, the gains become even larger when these representatives are chosen cleverly: for example, the maximum likelihood parameters of a DAG are generally much easier to compute than those of a graph which also includes bidirected edges, so by picking a DAG from each class that contains one, we can avoid many relatively expensive score computations on more complex graphs.

When we are looking for maximum likelihood parameters for an HTC-identifiable graph $G$ but our maximum likelihood fitting procedure has difficulty converging, it may be beneficial to apply the procedure to an `easier' \eqly{} equivalent graph $G'$ instead. After maximum likelihood parameters $(\Lambda', \Omega')$ have been found, we can compute $\Sigma = \phi_{G'}(\Lambda', \Omega')$, the point in $\cM(G')$ where the likelihood is maximized; this will be the same point for all \eqly{} equivalent models (up to the nongeneric case where $\Sigma \notin \cM(G)$---though then points arbitrarily close to $\Sigma$ will be included in $\cM(G)$).
For this $\Sigma$, we can compute parameters $(\Lambda, \Omega)$ for which $\phi_G(\Lambda, \Omega) = \Sigma$ using the algorithm in the proof of Theorem~1 of \citet{FoygelDraismaDrton2012_htc}. 

Because Theorem~\ref{thm:equiv_sufficient} only provides a sufficient condition for \eqic{} equivalence, an automatic procedure for model selection based on this theorem may sometimes fail to recognize that two classes of graphs are equivalent, and consider both separately (unless the results of Theorem~\ref{thm:allconstraints} are also considered by the algorithm, similar to how we used them above). While this means that some redundant computational work is done, it does not hurt the quality of model selection, and the gains in computation time are still enormous compared to testing all models individually: 34752 acyclic mixed graphs on four nodes are grouped into 419 clusters by Theorem~\ref{thm:equiv_sufficient}, which is very close to the 389 \eqic{} equivalence classes we would find by also looking at the constraints.

\subsection{FINITE-TO-ONE GRAPHS AND INEQUALITY CONSTRAINTS}\label{sec:ineq}

\begin{figure}[t]
  \centering
  \stackunder{\begin{tikzpicture}
    \node [circle,fill=black,inner sep=1pt] (a) at (0.0,1.0) [label=180.0:$\mathstrut a$] {};
    \node [circle,fill=black,inner sep=1pt] (b) at (1.0,2.0) [label=0.0:$\mathstrut b$] {};
    \node [circle,fill=black,inner sep=1pt] (c) at (1.0,1.0) [label=0.0:$\mathstrut c$] {};
    \node [circle,fill=black,inner sep=1pt] (d) at (1.0,0.0) [label=0.0:$\mathstrut d$] {};
    \draw [magenta!80!black,double,arrows=
      {Stealth[reversed,sep,length=.64ex,inset=.25ex,width=.8ex]-Stealth[sep,length=1ex]}]
      (a) -- (b);
    \draw [magenta!80!black,double,arrows=
      {Stealth[reversed,sep,length=.64ex,inset=.25ex,width=.8ex]-Stealth[sep,length=1ex]}]
      (a) -- (c);
    \draw [magenta!80!black,double,arrows=
      {Stealth[reversed,sep,length=.64ex,inset=.25ex,width=.8ex]-Stealth[sep,length=1ex]}]
      (a) -- (d);
  \end{tikzpicture}}{(a)}
  \quad
  \stackunder{\begin{tikzpicture}
    \node [circle,fill=black,inner sep=1pt] (a) at (0.0,1.0) [label=180.0:$\mathstrut a$] {};
    \node [circle,fill=black,inner sep=1pt] (b) at (1.0,2.0) [label=0.0:$\mathstrut b$] {};
    \node [circle,fill=black,inner sep=1pt] (c) at (1.0,1.0) [label=0.0:$\mathstrut c$] {};
    \node [circle,fill=black,inner sep=1pt] (d) at (1.0,0.0) [label=0.0:$\mathstrut d$] {};
    \draw [magenta!80!black,double,arrows=
      {Stealth[reversed,sep,length=.64ex,inset=.25ex,width=.8ex]-Stealth[sep,length=1ex]}]
      (a) -- (b);
    \draw [magenta!80!black,double,arrows=
      {Stealth[reversed,sep,length=.64ex,inset=.25ex,width=.8ex]-Stealth[sep,length=1ex]}]
      (a) -- (c);
    \draw [red,dashed,arrows=
      {Stealth[sep,length=1ex]-Stealth[sep,length=1ex]}]
      (a) -- (d);
    \draw [blue,arrows=
      {_-Stealth[sep,length=1ex]}]
      (c) -- (d);
  \end{tikzpicture}}{(b)}
  \quad
  \stackunder{\begin{tikzpicture}
    \node [circle,fill=black,inner sep=1pt] (a) at (0.0,1.0) [label=180.0:$\mathstrut a$] {};
    \node [circle,fill=black,inner sep=1pt] (b) at (1.0,2.0) [label=0.0:$\mathstrut b$] {};
    \node [circle,fill=black,inner sep=1pt] (c) at (1.0,1.0) [label=0.0:$\mathstrut c$] {};
    \node [circle,fill=black,inner sep=1pt] (d) at (1.0,0.0) [label=0.0:$\mathstrut d$] {};
    \draw [magenta!80!black,double,arrows=
      {Stealth[reversed,sep,length=.64ex,inset=.25ex,width=.8ex]-Stealth[sep,length=1ex]}]
      (a) -- (b);
    \draw [red,dashed,arrows=
      {Stealth[sep,length=1ex]-Stealth[sep,length=1ex]}]
      (a) -- (c);
    \draw [blue,arrows=
      {_-Stealth[sep,length=1ex]}]
      (b) -- (c);
    \draw [red,dashed,arrows=
      {Stealth[sep,length=1ex]-Stealth[sep,length=1ex]}]
      (a) -- (d);
    \draw [blue,arrows=
      {_-Stealth[sep,length=1ex]}]
      (c) -- (d);
  \end{tikzpicture}}{(c)}
  \quad
  \stackunder{\begin{tikzpicture}
    \node [circle,fill=black,inner sep=1pt] (a) at (0.0,1.0) [label=180.0:$\mathstrut a$] {};
    \node [circle,fill=black,inner sep=1pt] (b) at (1.0,2.0) [label=0.0:$\mathstrut b$] {};
    \node [circle,fill=black,inner sep=1pt] (c) at (1.0,1.0) [label=0.0:$\mathstrut c$] {};
    \node [circle,fill=black,inner sep=1pt] (d) at (1.0,0.0) [label=0.0:$\mathstrut d$] {};
    \draw [magenta!80!black,double,arrows=
      {Stealth[reversed,sep,length=.64ex,inset=.25ex,width=.8ex]-Stealth[sep,length=1ex]}]
      (a) -- (c);
    \draw [red,dashed,arrows=
      {Stealth[sep,length=1ex]-Stealth[sep,length=1ex]}]
      (a) -- (b);
    \draw [blue,arrows=
      {_-Stealth[sep,length=1ex]}]
      (c) -- (b);
    \draw [red,dashed,arrows=
      {Stealth[sep,length=1ex]-Stealth[sep,length=1ex]}]
      (a) -- (d);
    \draw [blue,arrows=
      {_-Stealth[sep,length=1ex]}]
      (c) -- (d);
  \end{tikzpicture}}{(d)}
\caption{Four acyclic finite-to-one graphs.}\label{fig:finitetoone}
\end{figure}
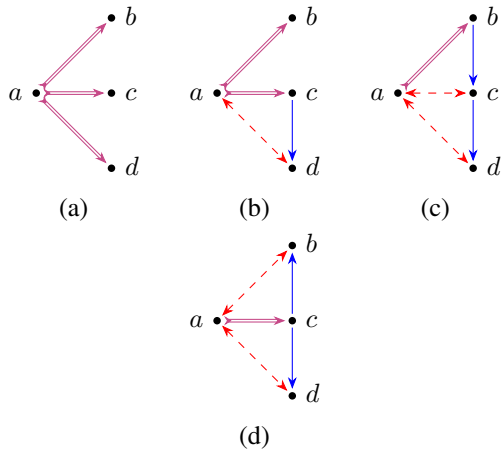

Among the acyclic graphs on four nodes, four (up to graph isomorphism) are HTC-inconclusive. These graphs are shown in Figure~\ref{fig:finitetoone}. By \citep[Table~1]{FoygelDraismaDrton2012_htc}, these graphs are generically finite-to-one but not identifiable. We used this in Section~\ref{sec:modelclasses} to assign them to \eqic{} equivalence classes, and found that all are \eqly{} equivalent to the saturated model.

Because these graphs are not generically identifiable, Theorem~\ref{thm:allconstraints} does not apply to them, so we cannot rule out the possibility that they impose inequality constraints. Indeed, the proposition below shows that one of the graphs imposes such a constraint, so that its model differs from the saturated model by a subset of nonzero measure. We believe the other graphs listed here impose similar constraints, and expect that the same is true for many other graphs that are generically finite-to-one but not generically identifiable.
\begin{proposition}\label{prop:ineq}
The model of the graph in Figure~\ref{fig:finitetoone}(a) imposes the inequality constraint $\rho_{bc.a} \cdot \rho_{cd.a} \cdot \rho_{bd.a} \leq 0$.
\end{proposition}


\section{EXPERIMENTAL RESULTS}\label{sec:experimental}

In Section~\ref{sec:experimental_Y}, we will describe the results of a model selection experiment using \eqic{} equivalence classes, choosing the class with the best BIC score from among those found in Section~\ref{sec:modelclasses}. We define the BIC score of an equivalence class as the best score among its models; this can be determined by computing the maximum likelihood of just one member model. RICF \citep{DrtonEichlerRichardson2009_RICF} was used to find maximum likelihood parameters. We will first describe some empirical results about its convergence behaviour.

\subsection{CONVERGENCE BEHAVIOUR OF RICF}

As already mentioned in Section~\ref{sec:consequences}, knowing the \eqic{} equivalence classes is very useful in a model selection problem. Here we elaborate on this, based on Monte Carlo results.

RICF finds maximum likelihood parameters in its first iteration if the graph is a DAG \citep{DrtonEichlerRichardson2009_RICF}. For all \eqic{} equivalence classes not containing a DAG (so all graphs in the class contain a bidirected edge), we saw evidence of local optima: when RICF was run with the same random data on different graphs in the same class, or with different initialization values, there would be different runs that all reported convergence but achieved different likelihoods. This could be addressed by random restarts. However, \citet{DrtonRichardson2004_LikelihoodMultimodality} observe that for one of these graphs, local optima only present themselves when the model is misspecified. If the same is true for all graphs with bidirected edges, then it would follow that local minima will not change the results of model selection, as they only affect the scores of models that would not have scored well anyway. Either way, it is clearly advantageous to run RICF on a DAG, for every \eqic{} equivalence class that contains one.

For equivalence classes where all graphs contain a bow, RICF failed to converge relatively often (in fact, \citet{DrtonEichlerRichardson2009_RICF} only discuss bow-free graphs). These equivalence classes contain graphs having different skeletons; as an example, the nine graphs represented by the pattern in Figure~\ref{fig:extra_equiv}(a) have a different skeleton than those in Figure~\ref{fig:extra_equiv}(b), but all are \eqly{} equivalent to each other. We found that for many data sets,
RICF failed to converge on all graphs in such a class having one skeleton, even with random restarts, but did converge easily on graphs with a different skeleton. Knowing the \eqic{} equivalence classes is a great benefit here, as it would be very difficult to determine scores for all graphs in this situation otherwise. When scoring a model, if we find that RICF fails to converge on a graph of this type (e.g., one from Figure~\ref{fig:extra_equiv}(a)), then we run RICF again on an \eqly{} equivalent graph with a different skeleton (one from Figure~\ref{fig:extra_equiv}{b}).


\subsection{USING \EQIC{} EQUIVALENCE CLASS SELECTION TO DETECT Y-STRUCTURES}\label{sec:experimental_Y}

\begin{figure}[t]
  \centering
  \stackunder{\begin{tikzpicture} 
    \node [circle,fill=black,inner sep=1pt] (a) at (0.0,0.866025403784) [label=180.0:$\mathstrut a$] {};
    \node [circle,fill=black,inner sep=1pt] (b) at (1.0,0.866025403784) [label=0.0:$\mathstrut b$] {};
    \node [circle,fill=black,inner sep=1pt] (c) at (0.5,0.0) [label=0.0:$\mathstrut c$] {};
    \node [circle,fill=black,inner sep=1pt] (d) at (0.5,-1.0) [label=0.0:$\mathstrut d$] {};
    \draw [blue,arrows=
      {_-Stealth[sep,length=1ex]}]
      (c) -- (d);
    \draw [green!67!black,arrows=
      {Stealth[sep,length=1ex]-Bar[sep,width=1ex] __}]
      (c) -- (a);
    \draw [green!67!black,arrows=
      {Stealth[sep,length=1ex]-Bar[sep,width=1ex] __}]
      (c) -- (b);
  \end{tikzpicture}}{(a)}
  \quad
  \stackunder{\begin{tikzpicture} 
    \node [circle,fill=black,inner sep=1pt] (a) at (0.0,0.866025403784) [label=180.0:$\mathstrut a$] {};
    \node [circle,fill=black,inner sep=1pt] (b) at (1.0,0.866025403784) [label=0.0:$\mathstrut b$] {};
    \node [circle,fill=black,inner sep=1pt] (c) at (0.5,0.0) [label=0.0:$\mathstrut c$] {};
    \node [circle,fill=black,inner sep=1pt] (d) at (0.5,-1.0) [label=0.0:$\mathstrut d$] {};
    \draw [blue,arrows=
      {_-Stealth[sep,length=1ex]}]
      (c) -- (d);
    \draw [green!67!black,arrows=
      {Stealth[sep,length=1ex]-Bar[sep,width=1ex] __}]
      (c) -- (b);
    \draw [red,dashed,arrows=
      {Stealth[sep,length=1ex]-Stealth[sep,length=1ex]}]
      (a) -- (c);
    \draw [red,dashed,arrows=
      {Stealth[sep,length=1ex]-Stealth[sep,length=1ex]}]
      (a)  .. controls +($.5*(0, -1.0)$) and +($.5*(-.5, 0.866025403784)$) .. (d);
  \end{tikzpicture}}{(b)}
\caption{(a) A Y-structure; (b) an extended Y-structure.}\label{fig:Ystruc}
\end{figure}
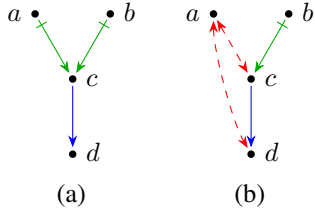

To demonstrate the practical usefulness of \eqic{} equivalence class selection, we consider the setup of \citet{MooijCremers_UAI2015CI}.\footnote{The code for reproducing these results is available online at \url{https://github.com/caus-am/aelsem}.} In a simulated dataset with $p \in \set{10, 30, 50}$ variables (with acyclic ground truth), they look at each ordered 4-tuple of distinct nodes, and use several independence tests to detect whether or not these nodes form a \emph{Y-structure} or an \emph{extended Y-structure} (see Figure~\ref{fig:Ystruc}). These are two of the simplest Markov equivalence classes (coinciding with \eqic{} equivalence classes) that must contain a directed edge in a fixed place that is not part of a bow, so detecting these structures in observational data allows us to draw conclusions about the results of interventions.

\citeauthor{MooijCremers_UAI2015CI} observed that detecting 4-tuples which were in either of these two classes yielded poor precision, especially for larger $p$. Precision improved when additional tests were added so that only Y-structures were detected.

\begin{figure}[t]
  \centering
  \includegraphics[width=.98\linewidth]{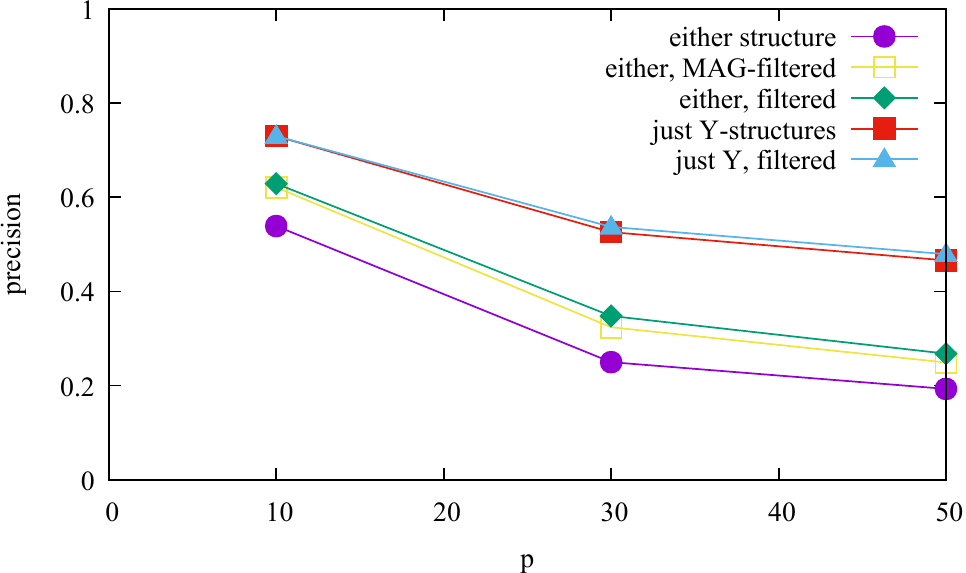}
  \caption{Precision of (extended) Y-structure detection.}\label{fig:precision}
\end{figure}

We applied \eqic{} equivalence class selection to this problem by taking the 4-tuples that tested positively according to one of these sets of tests, and then filtering out those for which the \eqic{} equivalence class with the best BIC score was different from the class / pair of classes being tested for. The resulting precisions are shown in Figure~\ref{fig:precision}.
The gains of this filtering procedure are significant 
when testing for both classes together: the precision is close to $0.1$ larger for all $p$. We also show the results of filtering using only Maximal Ancestral Graphs (MAGs) \citep{RichardsonSpirtes2002_MAGs}, thus only computing BIC scores of classes which can be described using vanishing partial correlations. This already yields a large improvement, but the gains from considering all \eqic{} equivalence classes instead are still significant, especially for the larger $p$. 
On the other hand, for 4-tuples that passed the more stringent Y-structure tests, the benefit of filtering is much smaller. Importantly, these gains in precision came at a very small cost in recall: of the true positives detected by a set of independence tests, at least 98\% passed through the filter for each $p$.

\section{CONCLUSION}\label{sec:conclusion}

We proposed the notion of \eqic{} equivalence for use in model selection among linear structural equation models, and showed how graphical criteria can be used to efficiently establish the equivalence of many models. Our experimental results show that the finer granularity of the resulting equivalence classes sometimes allows for improvements in model selection accuracy, compared to methods that only consider vanishing partial correlations.

While the experimental model selection results we show are limited to four nodes, we believe that the theoretical results described here can form the foundation of efficient causal inference algorithms on (much) larger numbers of nodes. For this, more work will need to be put into finding an efficient algorithm, for example by using ideas from \citet{Chickering2002_GES}.

\subsubsection*{Acknowledgements}

This project has received funding from the European Research Council (ERC) under the European Union's Horizon 2020 research and innovation programme
(grant agreement n\textordmasculine{} 639466).
The authors would like to thank Patrick Forr\'e, Tom Claassen and the anonymous reviewers for their valuable comments.

\subsubsection*{References}

\bibliographystyle{plainnat}
\DeclareRobustCommand{\VAN}[3]{#3 #1}
\begingroup 
\renewcommand{\section}[2]{}%

\endgroup


\clearpage
\appendix

\section*{SUPPLEMENTARY MATERIAL}

\section{PROOFS}\label{app:proofs}

\begin{proof}[Proof of Theorem~\ref{thm:allconstraints}]
  First assume $\Sigma \in \cM(G)$, and consider the equations
  \begin{equation}\label{eq:htc9.2}
    [(I - \Lambda)\trans \Sigma (I - \Lambda)]_{vw} = 0
    \quad 
    \text{for all $\set{v, w} \notin B$}
  \end{equation}
  (this is equation (9.2) in \citet{FoygelDraismaDrton2012_htc}).
  Treating the entries of $\Lambda$ and $\Sigma$ as variables $\lambda$ and $\sigma$, define the rational function $\delta \isdef \det(I - \Lambda)^{-1}$. Then the left-hand sides of \eqref{htc9.2} are polynomials in $\R[\lambda,\sigma,\delta]$ generating an ideal $\cI$.
  %
  This ideal consists of all polynomials that are zero everywhere on the set
  of all $(\lambda, \sigma)$ with $\delta$ well-defined and $\Sigma = \phi_G(\Lambda, \Omega)$ for some $\Omega \in \PD(B)$
  \citep[Section~8 of the supplement]{FoygelDraismaDrton2012_htc}.
  In addition to the left-hand sides of \eqref{htc9.2}, (using that $G$ is rationally identifiable) the ideal $\cI$ contains for each directed edge $(u,v) \in D$ an element $a(\sigma) \lambda_{uv} - b(\sigma)$ for some nonzero polynomials $a, b \in \R[\sigma]$. By multiplying a left-hand side of \eqref{htc9.2} by some polynomial and adding polynomial multiples of the ideal elements mentioned above, we can eliminate $\lambda$ and find an element in $\cI \cap \R[\sigma]$. This represents a polynomial constraint that \emph{all} $\Sigma \in \cM(G)$ must satisfy, and we see
  that the polynomial constraints described in the main text can be obtained in this way.
  %
  This implies that also the rational constraints \eqref{allconstraints} are satisfied for any $\Sigma$ that is additionally in the domain of $\Lambda_\cY$.

  If $\Sigma \in \cM(G)$, then there exist $\Lambda \in \R^D_{\text{reg}}$ and $\Omega \in \PD(B)$ for which $\Sigma = \phi_G(\Lambda, \Omega)$. In the generic case, $\Lambda_\cY$ will recover the parameter matrix $\Lambda$; in particular, $\Lambda_\cY(\Sigma) \in \R^D_{\text{reg}}.$


  Conversely, suppose $\Sigma$ is in the domain of $\Lambda_\cY$, $\Lambda_\cY(\Sigma) \in \R^D_{\text{reg}}$, and $\Sigma$ satisfies \eqref{allconstraints}. Let $\Lambda = \Lambda_\cY(\Sigma)$ and $\Omega = (I - \Lambda)\trans \Sigma (I - \Lambda)$. $\Sigma$ and $\Lambda$ again satisfy \eqref{htc9.2}:

  \emph{For $\set{v,w}$ with $v \notin Y_w$ and $w \notin Y_v$:} by \eqref{allconstraints};

  \emph{For $\set{v,w}$ with $w = y_i \in Y_v \cap \htr(v)$:}
    \begin{multline*}
      \left[ (I - \Lambda)\trans \Sigma \Lambda \right]_{wv}
      = (\bA \cdot \Lambda_{\pa(v),v})_i
      = \bb_i\\
      = \left[ (I - \Lambda)\trans \Sigma \right]_{wv}
    \end{multline*}
    (with $\bA$ and $\bb$ as in \citet{FoygelDraismaDrton2012_htc} / in Section~\ref{sec:prelim_htc}), which implies \eqref{htc9.2} for $\set{v,w}$;

    \emph{For $\set{v,w}$ with $w \in Y_v \setminus \htr(v)$:} Fix any node $v$, and let $W = Y_v \setminus \htr(v)$ and $I = \set{i \given y_i \in W}$. Then
    \begin{equation*}
      \left[ \Sigma \Lambda \right]_{Wv}
      = (\bA \cdot \Lambda_{\pa(v),v})_I
      = \bb_I
      = \Sigma_{Wv},
    \end{equation*}
    from which follows
    \begin{align*}
      \mathbf{0} &= \left[\Sigma (I - \Lambda) \right]_{Wv}\\
      &= \left[(I - \Lambda)\invtrans (I - \Lambda)\trans \Sigma (I - \Lambda) \right]_{Wv}\\
      &= \left[(I - \Lambda)\invtrans \Omega \right]_{Wv}
      = \sum_{x \in V} [(I - \Lambda)\invtrans ]_{Wx} \Omega_{xv}.
    \end{align*}
    In the case distinction that follows, we use that for all $x \neq v$ with $x \notin \sib(v)$, we have $x \in Y_v$, $v \in Y_x$, or neither.
    For $x$ in one of the sets $Y_v \cap \htr(v)$, $\set{x \given v \in Y_x \cap \htr(x)}$, or $\set{x \given x \neq v, x \notin \sib(v), x \notin Y_v, v \notin Y_x}$, we know from the previous two cases that $\Omega_{xv} = \Omega_{vx} = 0$. For $x = v$ or $x \in \sib(v)$, $[(I - \Lambda)\invtrans  ]_{wx} = 0$ for any $w \in W$, because such a directed path from such $x$ to $w$ would form a half-trek from $v$ to $w$. Similarly, for $x$ such that $v \in Y_x \setminus \htr(x)$, also $[(I - \Lambda)\invtrans]_{wx} = 0$, or there would again be a half-trek from $v$ to $w$.
    The above covers all cases except $x \in W$.
    So all terms in the previous sum with $x \notin W$ are zero, and we get
    \begin{equation*}
      \mathbf{0} = [(I - \Lambda)\invtrans ]_{WW} \Omega_{Wv}.
    \end{equation*}
    The entries of $[(I - \Lambda)\invtrans ]_{WW}$ are polynomials in the variables $\Lambda_{ij}$ and $\delta \isdef \det(I - \Lambda)^{-1}$.
    The matrix $(I - \Lambda)^{-1}$ equals $\delta$ times the adjugate matrix of $(I - \Lambda)$; this adjugate matrix includes the constant term $1$ in each diagonal entry, and no constant terms in any other entries.
    We see that the determinant of $[(I - \Lambda)\invtrans ]_{WW}$
    contains the term $\delta^{\lvert W \rvert}$ so it is not the zero polynomial. 
    It follows that $[(I - \Lambda)\invtrans ]_{WW}$ is generically invertible, 
    so the only solution to the above equation is $\Omega_{wv} = 0$ for all $w \in W$. This again implies \eqref{htc9.2}.

    Then $\Omega = (I - \Lambda)\trans \Sigma (I - \Lambda)$ is positive definite (using that $\Sigma$ is positive definite and $I - \Lambda$ is invertible) 
    and obeys $\Omega_{vw} = 0$ for $\set{v,w} \notin B$. For these parameters, $\Sigma = \phi_G(\Lambda, \Omega)$, showing that $\Sigma \in \cM(G)$.

  For generic $\Sigma \in \cM(G)$, the matrices $\bA$ that occur in the HTC-identifiability algorithm are invertible by \citep[Lemma~2]{FoygelDraismaDrton2012_htc}, so the assumption above that $\Sigma$ is in the domain of $\Lambda_\cY$ holds in the generic case. 
\end{proof}

\begin{proof}[Example of a constraint that is not HT-overidentifying]
  Figure~\ref{fig:overid_incomplete} 
\begin{figure}[b]
  \centering
  \begin{tikzpicture}
    \node [circle,fill=black,inner sep=1pt] (a) at (0.0,0.0) [label=270.0:$\mathstrut a$] {};
    \node [circle,fill=black,inner sep=1pt] (b) at (1.0,0.0) [label=270.0:$\mathstrut b$] {};
    \node [circle,fill=black,inner sep=1pt] (c) at (2.0,0.0) [label=270.0:$\mathstrut c$] {};
    \node [circle,fill=black,inner sep=1pt] (d) at (3.0,0.0) [label=270.0:$\mathstrut d$] {};
    \draw [magenta!80!black,double,arrows=
      {Stealth[reversed,sep,length=.64ex,inset=.25ex,width=.8ex]-Stealth[sep,length=1ex]}]
      (b) -- (a);
    \draw [red,dashed,arrows=
      {Stealth[sep,length=1ex]-Stealth[sep,length=1ex]}]
      (b) -- (c);
    \draw [magenta!80!black,double,arrows=
      {Stealth[reversed,sep,length=.64ex,inset=.25ex,width=.8ex]-Stealth[sep,length=1ex]}]
      (c) -- (d);
  \end{tikzpicture}
\caption{Graph imposing an equality constraint that is not HT-overidentifying.}\label{fig:overid_incomplete}
\end{figure}
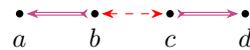
  gives an example graph which imposes one equality constraint, namely a vanishing tetrad constraint, which is not reported by Algorithm~2 of \citet{ChenTianPearl2014_LSEMOveridentify}.
  The graph is HTC-identifiable, so Theorem~\ref{thm:allconstraints} will report the constraint.
\end{proof}

\begin{proof}[Proof of Theorem~\ref{thm:equiv_sufficient}]
  Define the \emph{rank deficiency} of a graph as the difference between the number of columns and the generic rank of the Jacobian in \citep[proof of Theorem~2]{FoygelDraismaDrton2012_htc}. The rank deficiency of a graph is zero iff 
  that graph is finite-to-one, so the rank deficiency of $G'$ is zero. Removing a directed edge from a graph removes a column from the Jacobian, which will not increase the rank, so this reduces the rank deficiency by at most one; removing a bidirected edge adds a row, which will increase the largest order of a non-vanishing minor (this order equals the rank) by at most one. It follows that the rank deficiency of $G$ is at most $k$. In fact, it must be equal to $k$: suppose the rank deficiency of $G$ is $k' < k$. Then by deleting $k'$ columns from the Jacobian, a matrix of full column rank can be obtained. By deleting the corresponding $k'$ directed edges from $G$, a graph with rank deficiency 0 is obtained; however, this is a contradiction with the conditions of the theorem, which imply that all such graphs are infinite-to-one.

  By the statements in the beginning and end of the proof of Lemma~2 in the supplement of \citep{FoygelDraismaDrton2012_htc}, the dimensions of the images of $\phi_G$ and $\phi_{G'}$ (i.e.~of $\cM(G)$ and $\cM(G')$) equal the ranks of the respective Jacobians as considered above plus the respective numbers of (co)variance parameters in $\Omega$. It follows that the models' dimensions are equal.
  Because the models are defined by rational parameterizations, their Zariski closures are irreducible affine varieties \citep[Proposition~4.5.6]{CoxLittleOShea2015}. For two such varieties, $W \subsetneq V$ would imply that the dimension of $W$ is smaller than that of $V$.
  Because $\cM(G') \subseteq \cM(G)$, the models' Zariski closures must be equal; in other words, the models are \eqly{} equivalent.

  If additionally $G'$ imposes no inequality constraints, then $\cM(G')$ is equal to its Zariski closure up to a measure zero subset. Using again that $\cM(G') \subseteq \cM(G)$ and that $\cM(G)$ is contained in the models' common Zariski closure, we find that $\cM(G')$ and $\cM(G)$ are equal up to a measure zero subset.
\end{proof}

\begin{proof}[Proof of Proposition~\ref{prop:nowzohour}]
  First consider the special case of two acyclic bow-free graphs $G'_1, G'_2$ that differ in only one arrowhead which does not form a collider. 
  Then the graph $G$ containing the union of their edges (so having a bow where $G'_1$ and $G'_2$ had a different type of edge) is HTC-nonidentifiable, as the directed edge of the bow cannot be identified. An application of Corollary~\ref{cor:fin_equiv} shows that $G'_1$ and $G'_2$ are \eqly{} equivalent. 

  The \eqic{} equivalence of two arbitrary graphs $G'_1, G'_2$ with the same skeleton and colliders can now be shown using a finite sequence of such steps. Let $G^\Delta_1$ be the subgraph consisting of the edges of $G'_1$ that are not in $G'_2$, and let $W \subseteq V$ be the set of nodes where $G'_1$ and $G'_2$ have different incoming arrowheads. In $G^\Delta_1$, at each node $w \in W$,
  there is at most one arrowhead (otherwise there is a collider in either $G'_1$ or $G'_2$ that is not in the other)
  and at most one arrowtail (two tails would become two colliding heads in $G'_2$),
  so each node $w \in W$ has degree at most 2 in $G^\Delta_1$.

  Now find the finest partition of the edges in $G^\Delta_1$ that groups two edges together if they are incident at a node $w \in W$. The edge sets in this partition form paths that may only meet each other at the endpoints, at nodes not in $W$. None of these paths can have two directed edges pointing away from each other (with possibly other edges in between), as that would imply $G'_2$ contains a collider not in $G'_1$. So each of these paths is a half-trek. At each step, modify the first edge of such a half-trek by adding an arrowhead that appears in $G'_2$ or removing one that does not.
\end{proof}

\begin{proof}[Proof of Proposition~\ref{prop:ineq}]
  Assume w.l.o.g.~that $\Sigma$ is normalized to a correlation matrix (with ones on the diagonal). 
  Let $W = \set{b,c,d}$. Then for each pair of distinct nodes $x, y \in W$,
  \begin{align*}
    \sigma_{xy} &= \lambda_{ax}\lambda_{ay} + \lambda_{ax}\omega_{ay} + \omega_{ax}\lambda_{ay}\\
    &= (\lambda_{ax} + \omega_{ax})(\lambda_{ay} + \omega_{ay}) - \omega_{ax}\omega_{ay}\\
    &= \sigma_{ax}\sigma_{ay} - \omega_{ax}\omega_{ay}.
  \end{align*}
  It follows that
  \begin{equation}\label{eq:pair}
    \omega_{ax}\omega_{ay} = \sigma_{ax}\sigma_{ay} - \sigma_{xy}.
  \end{equation}
  We have three of these equations for the three possible pairs in $W$. If all three $\omega_{aw}$'s have the same sign, all these expressions will be positive; otherwise, two will be negative and one positive. Other combinations of positive and negative cannot be attained by any choice of $\Omega$. Using that \eqref{pair} has the same sign as $-\rho_{xy.a}$, we can write this as the inequality constraint $\rho_{bc.a} \cdot \rho_{cd.a} \cdot \rho_{bd.a} \leq 0$.

  The observed covariance matrix
  \begin{equation*}
    \Sigma =
    \begin{bmatrix}
      1 & 0 & 0 & 0 \\
      0 & 1 & 1/2 & 1/2 \\
      0 & 1/2 & 1 & 1/2 \\
      0 & 1/2 & 1/2 & 1
    \end{bmatrix}
  \end{equation*}
  is an example of a positive definite matrix (so an element of the saturated model) that does not satisfy this inequality constraint. In fact, it is easy to see that an open ball around $\Sigma$ has this property, showing that the inequality constraint rules out a part of the saturated model having nonzero measure.
\end{proof}

\clearpage
\section{ACYCLIC \EQIC{} EQUIVALENCE CLASSES ON FOUR NODES}\label{app:modelclasses}

Before presenting in Section~\ref{sec:modelclassestable} a full description of the \eqic{} equivalence classes of four-node acyclic graphs, we first introduce the additional graph notation we will use.

\subsection{NOTATION FOR SETS OF GRAPHS, CONTINUED}\label{sec:graph_notation_long}

In Section~\ref{sec:graph_notation}, we described how our graph patterns represent directed edges, bidirected edges, and bows.
The full list of markings we use in Section~\ref{sec:modelclassestable} is displayed in Figure~\ref{fig:legend_long}.
To describe these markings, we need to extend the well-known concept of colliders (see e.g.~\citet{SpirtesGlymourScheines2000}) 
to graph skeletons. A 2-edge path in the skeleton $S(G)$ of $G$ along three distinct nodes $(v_1, v_2, v_3)$ is called a \emph{collider} if for both edges of the path, a corresponding edge exists in the original graph $G$ with an arrowhead at $v_2$. Such an edge could be either bidirected, or directed towards $v_2$. We distinguish
two cases: if \emph{all} corresponding edges in the original graph $G$ have arrowheads at $v_2$ (in other words, there are no directed edges in $G$ from $v_2$ to either $v_1$ or $v_3$), we call the path a \emph{full collider}; otherwise we call it a \emph{partial collider}.
An example of a partial collider is the path $(a,b,c)$ in the instrumental variable model (Figure~\ref{fig:instrumental}(a)).
When looking at the endpoint at $v_2$ of an edge in $S(G)$ between $v_1$ and $v_2$, we say that it \emph{forms a (full/partial) collider} if a node $v_3$ exists for which the path along $(v_1, v_2, v_3)$ is a (full/partial) collider. Note that an endpoint can simultaneously form a full collider and a partial collider, for different choices of $v_3$.

In a graph pattern denoting a set of graphs, if all graphs in the set have the same edge(s) between two given nodes, this edge is displayed as for a single graph, as in Section~\ref{sec:graph_notation}. If the graphs in the set have different edges between two nodes, these are represented in the graph pattern by a green edge, which may have different endpoint markings. First of all, if all edges to be represented have an arrowhead at an endpoint, then the graph pattern will also show an arrowhead. Otherwise, the following special markings are used (illustrated in Figure~\ref{fig:legend_long}):
\begin{figure}[t]
  \centering
  \begin{tikzpicture} 
    \node (step) at (0, -.45) {};
    \draw ($(-10pt,0)-(step)$) rectangle ($(6.5,0)+7*(step)$);
    \node [circle,fill=black,inner sep=1pt] (a1) at (0,0) {};
    \node [circle,fill=black,inner sep=1pt] (b1) at (1,0) {};
    \node [circle,fill=black,inner sep=1pt] (a2) at ($(0,0)+(step)$) {};
    \node [circle,fill=black,inner sep=1pt] (b2) at ($(1,0)+(step)$) {};
    \node [circle,fill=black,inner sep=1pt] (a3) at ($(0,0)+2*(step)$) {};
    \node [circle,fill=black,inner sep=1pt] (b3) at ($(1,0)+2*(step)$) {};
    \node [circle,fill=black,inner sep=1pt] (a4) at ($(0,0)+3*(step)$) {};
    \node (b4) at ($(.6,0)+3*(step)$) {};
    \node [circle,fill=black,inner sep=1pt] (a5) at ($(0,0)+4*(step)$) {};
    \node (b5) at ($(.6,0)+4*(step)$) {};
    \node [circle,fill=black,inner sep=1pt] (a6) at ($(0,0)+5*(step)$) {};
    \node (b6) at ($(.6,0)+5*(step)$) {};
    \node (a7) at ($(0,0)+6*(step)$) {};
    \node (b7) at ($(1,0)+6*(step)$) {};
    \draw [blue,arrows=
      {Stealth[sep,length=1ex]-_}]
      (b1) -- (a1);
    \node [anchor=west] at (1.2,0) {directed edge};
    \draw [red,dashed,arrows=
      {Stealth[sep,length=1ex]-Stealth[sep,length=1ex]}]
      (a2) -- (b2);
    \node [anchor=west] at ($(1.2,0)+(step)$) {bidirected edge};
    \draw [magenta!80!black,double,arrows=
      {Stealth[reversed,sep,length=.64ex,inset=.25ex,width=.8ex]-Stealth[sep,length=1ex]}]
      (a3) -- (b3);
    \node [anchor=west] at ($(1.2,0)+2*(step)$) {bow: directed and bidirected edge};
    \draw [green!67!black,arrows=
      {-Bar[sep,width=1ex] __}]
      (b4) -- (a4);
    \node [anchor=west] at ($(1.2,0)+3*(step)$) {endpoint not forming collider};
    \draw [green!67!black,arrows=
      {-Bracket[reversed,sep,length=.4ex,width=.8ex]}]
      (b5) -- (a5);
    \node [anchor=west] at ($(1.2,0)+4*(step)$) {endpoint forming partial collider};
    \draw [green!67!black,arrows=
      {-Rays[n=6,sep,length=.9ex]}]
      (b6) -- (a6);
    \node [anchor=west] at ($(1.2,0)+5*(step)$) {other endpoint};
    \draw [gray,dotted]
      (a7) -- (b7);
    \node [anchor=west] at ($(1.2,0)+6*(step)$) {edge not present in all graphs};
  \end{tikzpicture}
\caption{Legend for the edges and endpoints we use to denote (sets of) graphs.}\label{fig:legend_long}
\end{figure}
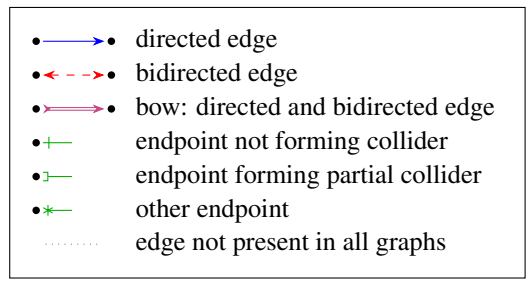
a plus sign if the endpoint never forms a collider in any of the graphs; a bracket if it always forms a partial collider; and a star in all other cases.
Finally, if two nodes are adjacent in some graphs in the set but not in others, then an edge is shown in dotted grey; other aspects of the edge's appearance are determined using only the graphs in which the nodes are adjacent.

\subsection{TABLE OF EQUIVALENCE CLASSES}\label{sec:modelclassestable}

Each entry in the table below represents one \eqic{} equivalence class and its isomorphisms (obtained by relabelling the nodes). The leftmost column shows the set of graphs in this equivalence class, using the notation for graph patterns explained in Sections~\ref{sec:graph_notation} and~\ref{sec:graph_notation_long}. For each equivalence class, its algebraic constraints (found by Theorem~\ref{thm:allconstraints} and converted to polynomial form) are also listed.

For many equivalence classes, the set of graphs that fit the pattern in the leftmost column is a superset of the equivalence class: some acyclic graphs that match the pattern do not actually belong in the class, but our notation is not expressive enough to show the precise inclusions in one figure. For each equivalence class where this is the case, a second column of smaller graph patterns appears. The union of the sets of graphs described by these patterns equals the equivalence class exactly. Because many of the graph patterns in the union are graph isomorphisms of each other, we show only one member of each isomorphism class, and indicate how many isomorphisms there are by writing `$n \times$'. (Each of these isomorphisms is an automorphism of the leftmost graph pattern, and there are always exactly $n$ such isomorphic patterns, so no information is lost by not listing them explicitly.) For example, the first equivalence class listed in Section~\ref{sec:table_dim3} consists of all acyclic graphs with at least one edge of any type between each pair of nodes in $\set{a,b,c}$, as well as all graphs in which one of those pairs has no edge between it, but the other two pairs make up a partial collider.

The equivalence classes are split out below according to their dimension. This refers to the dimension of the set of \emph{correlation} rather than covariance matrices, thus normalizing out the variance $\Sigma_{vv}$ associated with each node. This number has the convenient property that it equals the minimum number of edges among all graphs in an equivalence class.

\raggedbottom




\setlength{\tabulinesep}{3pt}

\setcounter{subsubsection}{-1}
\subsubsection{Dimension 0}

\begin{tabu}{cc}
  \hline

  \\*
  \multicolumn{2}{l}{$\sigma_{ab} = 0$, $\sigma_{ac} = 0$, $\sigma_{ad} = 0$,}
  \\*
  \multicolumn{2}{l}{$\sigma_{bc} = 0$, $\sigma_{bd} = 0$, $\sigma_{cd} = 0$}
  \\
  \hline
\end{tabu}

\subsubsection{Dimension 1}

\begin{tabu}{cc}
  \hline
  %
  \\*
  \multicolumn{2}{l}{$\sigma_{ac} = 0$, $\sigma_{ad} = 0$, $\sigma_{bc} = 0$, $\sigma_{bd} = 0$, $\sigma_{cd} = 0$}
  \\
  \hline
\end{tabu}

\subsubsection{Dimension 2}

\begin{tabu}{cc}
  \hline
  %
  \\*
  \multicolumn{2}{l}{$\sigma_{ac} = 0$, $\sigma_{ad} = 0$, $\sigma_{bd} = 0$, $\sigma_{cd} = 0$}
  \\
  \hline
\end{tabu}

\subsubsection{Dimension 3}\label{sec:table_dim3}

\begin{tabu}{cc}
  \hline
  \multirow{2}{*}{\parbox[t]{3cm}{
  %
  \\*
  \multicolumn{2}{l}{$\sigma_{bc} = 0$, $\sigma_{bd} = 0$, $\sigma_{cd} = 0$}
  \\
  \hline
\end{tabu}

\subsubsection{Dimension 4}

\begin{tabu}{cc}
  \hline
  %
  \\*
  \multicolumn{2}{l}{$\sigma_{cd} = 0$, $\sigma_{ac}\sigma_{bd} - \sigma_{ad}\sigma_{bc} = 0$}
  \\
  \hline
\end{tabu}

\subsubsection{Dimension 5}

\begin{tabu}{cc}
  \hline
  \multirow{2}{*}{\parbox[t]{3cm}{
  %
  \\*
  \multicolumn{2}{l}{$(p_1 \sigma_{aa} - p_2 \sigma_{ad})(\sigma_{cd}\sigma_{bb} - \sigma_{bc}\sigma_{bd})$}
  \\*
  \multicolumn{2}{l}{${} - (p_2 \sigma_{cd} - p_1 \sigma_{ac})(\sigma_{ab}\sigma_{bd} - \sigma_{ad}\sigma_{bb}) = 0$} 
  \\*
  \multicolumn{2}{l}{where $p_1 = \sigma_{bb}\sigma_{dd} - \sigma_{bd}^2$ and $p_2 = \sigma_{ad}\sigma_{bb} - \sigma_{bd}\sigma_{ab}$}
  \\
  \hline

  \\
  \hline
\end{tabu}

\subsubsection{Dimension 6}

\begin{tabu}{cc}
  \hline
  \multirow{2}{*}{\parbox[t]{3cm}{
  %
  \\
  \hline
\end{tabu}


\end{document}